\documentclass[oneside,english]{amsart}
\usepackage[latin1]{inputenc}
\usepackage{amssymb}

\makeatletter

\theoremstyle{plain}
 \newtheorem{thm}{Theorem}[section]
 \numberwithin{equation}{section} 
 \numberwithin{figure}{section} 
 \theoremstyle{remark}
 \newtheorem*{acknowledgement*}{Acknowledgement}
 \theoremstyle{plain}
 \newtheorem{conjecture}[thm]{Conjecture} 
 \theoremstyle{plain}
 \theoremstyle{plain}
 \newtheorem{lem}[thm]{Lemma} 
 \theoremstyle{definition}
  \newtheorem{question}[thm]{Question}
 \theoremstyle{plain}
 \newtheorem{prop}[thm]{Proposition} 
 \theoremstyle{plain}
 \newtheorem*{thm*}{Theorem}

\usepackage{color}
\usepackage{epsfig}
\newcommand{\hide}[1]{}

\newcommand{\ZZ}{\mathbb{Z}}
\newcommand{\NN}{\mathbb{N}}
\newcommand{\ZD}{{\mathbb{Z}^d}}
\newcommand{\GG}{\mathbb{G}}
\newcommand{\htop}{h_{\mathit{top}}}

\newcommand{\UDP}{unique derivation property}

\usepackage{babel}
\makeatother
\begin{document}

\title{Finite entropy for multidimensional cellular automata}
\author{Tom Meyerovitch}
\begin{abstract}
Let $X=S^\GG$ where $\GG$ is a countable group and $S$ is a finite
set. A cellular automaton (CA) is an endomorphism $T : X \to X$
(continuous, commuting with the action of $\GG$). Shereshevsky
\cite{shereshevsky_expansive_1993} proved that for $\GG=\ZD$ with
$d>1$ no CA can be forward expansive, raising the following
conjecture: For $G=\ZD$, $d>1$ the topological entropy of any CA is
either zero or infinite. Morris and Ward \cite{morris_ward_98},
proved this for linear CA's, leaving the original conjecture open.
We show that this conjecture is false, proving that for any $d$
there exist a $d$-dimensional CA with finite, nonzero topological
entropy. We also discuss a measure-theoretic counterpart of this
question for measure-preserving CA's.
\end{abstract}

\email{tomm@post.tau.ac.il}


\maketitle

\section{\label{sec:Introduction}Introduction}
Cellular automata form a class of dynamical systems which has been
extensively studied since the 1940's, going back to some work of von
Neumann \cite{vonneumann66} and others. For some survey papers on
different aspects of cellular automata see
\cite{BMS98,BM95,kari_survey_2005,wolfram02}. The following is an
intuitive description of a cellular automaton: Consider an infinite
mesh of finite state machines, called \emph{cells} interconnected
locally with each other. These cells change their states
synchronously depending on the states of some nearby cells, as
determined by a local update rule. All the cells use the same update
rule so that the system is homogeneous.

The \emph{entropy} of a cellular automaton is, informally, the rate
of information flow required to describe the long term evolution of
any finite number of cells. If the mesh of cells is a one
dimensional array (bi-infinite or one sided), a simple argument
shows that the entropy must be finite. Intuitively, when the cells
are arranged in a lattice of dimension $2$ or more, it can be
expected the entropy is either zero (the evolution of a finite
number of cells is eventually determined by the history of these
cells) or infinite (the rate of information to describe the
evolution of a set $F$ of cells grows to infinity with the
cardinality of $F$). This intuition, backed up by some supporting
results gave raise to the following conjecture
\cite{damico_2003,morris_ward_98}, attributed to Shereshevsky:
\begin{conjecture}
A $d$-dimensional cellular automaton, with $d > 1$, can not have
finite, positive entropy.
\end{conjecture}
\label{conj:trivial-CA-entropy} The main result of this paper is a
counter-example to this conjecture. For any $d \ge 2$ we describe a
$\ZD$ cellular automaton and give non-trivial bounds in its
topological entropy. An important tool implied in this construction
is a certain aperiodic set of tiles, associated to a substitution
system. The cellular automaton we describe was introduced by Kari
\cite{kari_xor_snake_1994} for $d=2$, to prove certain
undecidability results on cellular automata.

The paper is organized as follows: In section \ref{sec:Background}
we introduce notation and give brief definitions of cellular
automata, subshifts and entropy. Section \ref{sec:Background}
concludes with previous results and background on the entropy
problem for cellular automata. Section \ref{sec:CA-topo-entropy}
consists of a construction of an automaton with non-trivial
topological entropy, along with entropy estimations. Section
\ref{sec:Hilbet-path} contains an explanation of a key ingredient
for this construction - the so called ``Hilbert-Tiles''. Some tools
developed in sections \ref{sec:CA-topo-entropy} and
\ref{sec:Hilbet-path} are applied in section
\ref{sec:measure-entropy} to study measure-theoretic entropy of
surjective cellular automata. The last section of this paper
contains some further comments and open questions.

\textbf{Acknowledgements:} This work is part of the author's Ph.D.
carried under the supervision of Professor Jon Aaronson in Tel-Aviv
University. I would like to express my appreciation for his guidance
and encouragement. I would also like to thank T. Ward for helpful
suggestions and comments.

\section{\label{sec:Background}Preliminaries and background}

\subsection{Subshifts and cellular automata}
Let $S$ be some finite set, which we call the possible \emph{states}
of each cell. A \emph{configuration} is an element $x \in S^\GG$,
where $\GG$ is a countable set. For $n \in \GG$ We denote by $x_n$
the state of the cell at coordinate $n$ in the configuration $x$. If
$F \subset \GG$, we denote by $x_F$ the restriction of $x$ to $F$ -
the states of those cells located in $F$. We put the product
topology on $S^\GG$, which is compact and metrizable. If $\GG$ is a
semigroup, we can define an action $\sigma$ of the group $\GG$ on
$X$ by translations: For $h,g \in G$ and $x \in S^\GG$, $(\sigma_g
x)_{h}=x_{hg}$. Each map $\sigma_g$ is a continuous transformation
of $X$, and if $\GG$ is a group, $\{\sigma_g :~ g\in \GG\}$ are
homeomorphisms. The classic cases are $\GG=\ZZ$ and $\GG=\NN$, in
which case the action $\sigma$ is generated by one transformation,
called the \emph{shift map}. In this case, $S^\ZZ$ (or $S^\NN$) is
called a \emph{full shift} space (one sided or two sided,
respectively). If $\GG=\ZD$ where $d\ge 1$ then $S^\ZD$ is a
\emph{$d$-dimensional full shift}. A subset $X \subset S^\ZD$ which
is closed and invariant to the action of $\sigma$ is called a
$\ZD$-subshift. Any subshift can be described by a set of
\emph{forbidden blocks} $FB = \{ b_1, \ldots b_n, \ldots\}$ where
$b_i \in S^{F_i}$ and $F_i \subset \ZD$ are finite sets as follows:
$$X= \{x \in S^\ZD:~ x_{F_i+m} \ne b_i ~\forall m \in \ZD, i \in \NN \} $$
A subshift $X$ is called \emph{a subshift of finite type} (SFT) if
there exist a finite set $FB$ of forbidden blocks for $X$. We say
that a subshift $X \subset S^\ZD$ is a \emph{factor} of a subshift
$Y\subset \tilde{S}^\ZD$ (or $Y$ \emph{extends} $X$) if there exist
a surjective map $\pi: Y \to X$ which commutes with the shift action
$\sigma$.

For a countable group $\GG$, a \emph{$\GG$-cellular automaton} is a
pair $(S^\GG,T)$ where $T:S^\GG \to S^\GG$ is a continuous
transformation  which commutes with the $\GG$-action $\sigma$ of
translations. We abbreviate CA for cellular automaton. The
transformation $T$ is always of the form $(Tx)_n=t(x_{F+n})$ where
$F \subset \GG$ is a finite set, and $t:S^F \to S$ is called the
\emph{local rule} of the $CA$. Mostly, we restrict attention to the
case $\GG=\ZD$.

Often a wider definition of a $\ZD$-CA is  used, and any continuous,
shift commuting transformation of an SFT is called a $\ZD$-CA. In
this paper we deal specifically with CA on a full-shift.
\subsection{Topological and measure-theoretic entropy}
Let $(X,T)$ be a CA with $X=S^\GG$, and $F$ a finite subset of
$\GG$. We denote:
$$W(F,n,T):=\left\{\left( (x)_F,(Tx)_F,\ldots,(T^{n-1}x)_F \right):~ x \in
X \right\}$$ $W(F,n,T)$ is the set of possible configurations for
cells inside $F$ along $n$ iterations of $T$.

 The \emph{topological
entropy} of $(X,T)$ is:
$$\htop(X,T)=\sup_{F \subset \GG}\lim_{n \to
\infty}\frac{1}{n}\log|W(F,n,T)|$$ where the supremum above is over
all finite subsets $F$ of $\GG$, and the limit above exists by
subadditivity of the sequence $\log|W(F,n,T)|$.

 For $ w =(w_0,\ldots,w_{n-1}) \in W(F,n,T)$, let $[w]= \{x \in X :~ (T^ix)_F = w_i ~ 0 \le i < n\}$.
  Let $\mathcal{P}(X,T)$ be the set of probability measures on the Borel
$\sigma$-algebra $\mathcal{B}$ of $X$ which are invariant to $T$
($\mu(T^{-1}B)=\mu(B)$ any Borel $B \subset X$). For $\mu \in
\mathcal{P}(X,T)$. The measure-theoretic entropy of $(X,\mu,T)$ is:
$$h_{\mu}(X,\mathcal{B},T) =\sup_{F \subset \GG}\lim_{n \to
\infty}\frac{1}{n}H_n$$ with:
$$H_n = -\sum_{w \in W(F,n,T)}\mu([w])\log(\mu([w]))$$ where the supremum is again over all finite subsets $F$ of
$\GG$, and the limit above exists by subadditivity of the sequence
$H_n$. Both the quantities $\htop(X,T)$ and
$h_{\mu}(X,\mathcal{B},T)$ indicate the ''information flow'' of the
long term evolution of any finite number of cells by $T$. The
conceptual difference is that the topological entropy gives the rate
of information required to describe any evolution of finitely many
cells, where as $h_{\mu}(X,\mathcal{B},T)$ gives the rate of
information required to describe a ''typical'' evolution, where
''typical'' is with respect to the probability $\mu$.

An important result connecting the topological entropy and the
measure theoretic entropy, known as the \emph{variational principle}
states that: $$\htop(X,T)=\sup_{\mu \in
\mathcal{P}(X,T)}h_\mu(X,\mathcal{B},T)$$
 Entropy of a transformation (topological or measure theoretic) is defined in a
broader context. See \cite{MCK76} for definitions, proofs and a
detailed discussion in the context of compact spaces.

\subsection{Substitution systems and tilings}
A $\ZD$-\emph{tiling system} consists of a finite set $S$ of
''square tiles'' and some adjacency rules $R \subset S^F$ with
$F=\{0,\pm e_1,\ldots ,\pm e_d\}$, which determine when a tile $s_1$
is allowed to be placed next to a tile $s_2$ (and in which
directions). A configuration $x\in S^{F}$ for some $F \subset \ZD$
is \emph{valid} at $n \in F$ if the neighbors of the cell at $n$
obey the adjacency rules:  $x_{F+n} \in R$. Evidently, the set of
infinite valid configurations in $S^{\mathbb{Z}^2}$ are an SFT. We
call this SFT the associated subshift of $(S,R)$ and denote it by
$X_{R}$.

 A \emph{$\mathbb{Z}^d$-substitution system} consists
of a finite set $S$ and a substitution rule $\rho$ which is a
function $\rho:S \to S^{F_k}$ where $F_k=\{1,\ldots,k\}^d$. The
function $\rho$ naturally extends to a map $\rho^n:S \to
S^{F_{k^n}}$ by applying $n$ iterations of $\rho$, and to a map
$\rho: S^\ZD \to S^\ZD$. For $F \subset \ZD$ and a configuration
$x_F \in s^F$,  we say that $x_F$ is \emph{admissible} for the
substitution system $(S,\rho)$, if it appears as a sub-configuration
of $\rho^n(s)$ for some $s \in S$ and $n \in \NN$. To a substitution
system $(S,\rho)$ there is an associated subshift $X_\rho$, whose
forbidden blocks are the non-admissible blocks for $(S,\rho)$.

We say that a substitution system a has the \emph{\UDP} if any for
$x \in X_\rho$ there exist a unique $n \in F_k$ and $y \in X_\rho$
such that $x=\sigma_n\rho(y)$. Equivalently, this means that given a
part of a finite configuration $ x \in S^{F_{k^n}}$ such that
$x=\rho^n(s)$ for some $s\ \in S$ which has been shifted, the shift
can be recovered uniquely modulo $k\ZD$. Mozes \cite{mozes89}, based
on some earlier work of Robinson \cite{robinson71}, devised an
algorithm for implementing a subshift associated with a substitution
system via local constraints - a tiling system. Namely, Mozes proved
the following theorem:

\begin{thm}
\label{thm:subst}(Theorem 4.5 of \cite{mozes89}) Let $X_\rho \subset
S^\ZD$ be the $\ZD$-subshift associated to a substitution rule $s:S
\rightarrow S ^{F_{k}}$ with unique derivation. If $d>1$, then there
exists a tiling system $(\tilde{S},R)$ such that the associated
subshift $X_R$ extends $X_\rho$.
\end{thm}

To be precise, Mozes proved this theorem for $d=2$, but the proof
extends to any $d>1$. In fact, Goodman-Strauss \cite{GS98} proved a
more general result, relaying on Mozes' techniques, which implies
theorem \ref{thm:subst} for any $d>1$. For $d=1$, theorem
\ref{thm:subst} does not hold: there are classic examples of
one-dimensional substitution systems which are not factors of an
SFT.

\subsection{The CA entropy problem}
We turn to explain some motivation for conjecture
\ref{conj:trivial-CA-entropy}. A $\GG$-CA $(S^\GG,T)$ is
\emph{expansive} if there exist some finite  $F \subset \GG$ such
that for any $x \ne y \in S^{\GG}$ there exist $n \ge 0$ such that
$(T^nx)_F \ne (T^ny)_F$. Expansiveness of a CA means that by
observing the states of some finite number of cells under iterations
of $F$, we can eventually distinguish between any $2$
configurations. It can be shown that an expansive CA always has
finite, non-zero topological entropy. Shereshevsky
\cite{shereshevsky_expansive_1993} proved that a $\ZD$-CA with $d>1$
can not be expansive, which gave raise to conjecture
\ref{conj:trivial-CA-entropy}, stated in the introduction. Motivated
by this, Morris and Ward \cite{morris_ward_98} proved a result on
entropy of group automorphisms, which implies this conjecture for
some subclass of cellular automata called \emph{linear cellular
automata}. This result has been refined \cite{damico_2003}, to show
that linear CA
 are either ''sensitive to initial conditions'', in which case the
entropy is infinity, or equicontinuous, in which case the entropy is
zero. Lakshtanov and Langvagen \cite{LL2004} prove that any
multidimensional cellular automaton which admits a \emph{spaceship}
has infinite entropy. A \emph{spaceship} for a CA $(X,T)$ is
configuration $x \in X$ which differs from a $\sigma$-invariant
configuration in only finitely many cells, and such that $T^kx =
\sigma_{n} x$ for some $n \in \ZD \setminus \{0\}$ and some integer $k
>0$, but $x$ is not itself $\sigma$-invariant.

We point out that for any $d \ge 1$ there are trivial examples of continuous,
shift commuting transformations $T:X \to X$ with finite, positive entropy where $X \subset S^{\ZD}$ is an SFT:

Let $$X= \{ x \in S^{\ZD} :~ x_{n}=x_{n+e_i} ~\forall 1<i \le d\}$$
where $\{e_1,\ldots,e_d\}$ are the standard generators of $\ZD$
and $T(x)= \sigma_{e_1}x$ is the shift in the direction of $e_1$. In this case, the example is essentially $1$-dimensional, and the topological entropy of $T$ is $\log |S|$.

\section{\label{sec:CA-topo-entropy}Surjective $\ZD$ CA with finite nonzero entropy}
In this section we present a construction of a surjective $\ZD$-CA,
associated with a set of directed tiles $S$ and a finite group
$\Gamma$. This is a simple generalization of Kari's CA
\cite{kari_xor_snake_1994}. We prove that if $S$ has some special
properties, the associated CA has non-trivial entropy. The existence
of a set of directed tiles with the required properties is proved in
section \ref{sec:Hilbet-path}.

A set of \emph{directed tiles} is a tiling system with a forward
direction $d(s) \in \{\pm e_1,\ldots, \pm e_d\}$ associated to each
tiles $s \in S$. Given a configuration $x \in S^{\ZZ^2}$, a path
defined by $x$ is a sequence $p_1,p_2,\ldots$ with $p_n \in \ZD$
obtained by traversing the forward directions of $x$:
$p_{n+1}=p_n+d(x_{p_n})$. Given $x \in S^{\ZZ^2}$, a path $\{p_n\}$
is valid if $x$ is valid at every $p_n$. A set of directed tiles
$S$, such that no valid path in $x \in S^\ZD$ forms a loop is called
an \emph{acyclic} set of tiles.

Fix a finite group $\Gamma$. The group operation of $\Gamma$ is
written here in additive notation. Let $S$ be a set of directed set
of tiles. Recall that for $s \in S$, $d(s)$ denotes the forward
direction of the tile $s$.

 Define
$T_S:X \to X$ with $X=(S \times \Gamma)^{\ZD}$ as follows:

\begin{equation}
\label{eq:kari-ca}
T_{S,\Gamma}(x,y)_n:=\left\{ \begin{array}{cc} (x_n,y_n+y_{n+d(x_n)}) & \mbox{if $x$ is valid at $n$}\\
 (x_n,y_n) & \mbox{otherwise} \end{array} \right.
 \end{equation}
 In this cellular automaton cells ``transmit-information'' along
 valid paths. We prove that when  the number of infinite valid paths
 is bounded, the topological entropy is finite.

 For $x \in X$, denote by $\gamma(x)_n$ the $\Gamma$-part of state
of the cell at $n$ in $x$, and by $s(x)_n$ the $S$-part of the
state.
\begin{lem}
\label{lem:kari-ca-surjective} With $S$ an acyclic set of tiles as
above, and $T_{S,\Gamma}$ defined according to equation
\eqref{eq:kari-ca}, $T_{S,\Gamma}$ is surjective.
\end{lem}
\begin{proof}
It is sufficient to prove that the image of $T_{S,\Gamma}$ is dense
in $X$, that is, to any finite  $F \subset \ZD$ and any $y \in X$
there is some $x \in X$ such that $(T_{S,\Gamma}x)_F=y_F$. Let $y
\in X$ and $F \subset \ZD$ a finite set.  We describe  $x \in X$ as
follows:
 Set $s(x)_n=s(y)_n$ for any $n \in \ZD$.  For $n \in \ZD \setminus F$,
 set $\gamma(x)_n=0$ (the identity of the group $\Gamma$).
 It remains to define the $\gamma(x)_n$ for $n\in F$. We define $\gamma(x)_n$ iteratively, according to the length of maximal valid path in $x$ beginning at
 $n$, which remains in $F$. Since $S$ is acyclic and $F$ is finite, any such path must
 be finite.
 If $x$ is not valid at $n$, or the successor of $n$ in $x$ is not in $F$, define $\gamma(x)_n=\gamma(y)_n$.
 Suppose now that the length of the maximal valid path in $x$
 beginning at $n$ is $k>1$, and we have defined
 $\gamma(x)_{n+d(x_n)}$, So we define
 $\gamma(x)_n=\gamma(y)_n-\gamma(x)_{n+d(x_n)}$. It can now be
 verified that indeed $(T_{S,\Gamma}x)_F=y_F$.
\end{proof}

Suppose $S$ is an acyclic set of tiles, and $w \in S^\ZD$. Let $X_w
= \{ (x,y) \in S^\ZD\times\Gamma^\ZD:~ x=w\}$ and identify it with
naturally as a subset of $X$. The set $X_w$ is a closed subset of
$X$, and is invariant under $T_{S,\Gamma}$.  Let us define a
directed graph $G_w=(V_w,E_w)$ with vertex set $V_{w}=\ZD$ and the
edges
$$E_{w}=\{(n,n+d(w_n)):~ n \in \ZD,\, \mbox{$w$ is valid at $n$}\}$$
We say that $K \subset \ZD$ is connected in $G_w$ if for any two
cells in $K$ there is a directed path in $G_w$ from one to the other
(but not necessarily in both directions). A connected component of
$G_w$ is a connected set in $G_w$ which is maximal with respect to
inclusion. For any connected component $K$ of $G_w$, and any $x \in
X_w$, the CA $T_{S,\Gamma}$ acts on $x_K$ independently of the
states outside $K$.

Fix $w \in S^{\ZD}$. Suppose $G_w$ has $c$ infinite components,
denoted by $K_1,\ldots,K_c$. Let $K_0=\ZD \setminus
\bigcup_{i=1}^{c}K_c$ denote the union of cells which are not part
of forward infinite valid path in $w$. Let $(X_{i},T_i)$ denote the
system corresponding to those cells in $K_i$ with the action of
$T_{S,\Gamma}$. The systems $(X_{i},T_i)$ are factors of
$(X_w,T_{S,\Gamma})$.

 For $w \in S^{\ZD}$, denote by $I(w)$ the number of infinite connected components of $G_w$, which is equal to
the maximal number of pairwise disjoint, forward infinite valid
paths in $w$. For a directed set of tiles $S$, let:
$$I(S)=\sup_{w \in S^\ZD}I(w)$$

\begin{lem}
\label{lem:no-inifite-path-zero-entropy} Let $(X_0,T_0)$  be the
system corresponding to the cells which are not part of forward
infinite valid paths as above, $\htop(X_0,T_0)=0$.
\end{lem}
\begin{proof}
Since any valid path in $X_0$ is forward finite, a simple induction
on the length of the path shows that the sates of the cells in a
valid path of length at most $k$ is $2^k$-periodic. Thus,
$(X_0,T_{0})$ is isomorphic to an inverse limit of (finite) periodic
systems, and thus has $0$ topological entropy.
\end{proof}

\begin{lem}
\label{lem:infinite-path-entropy} Let $X_w$,and $T_{S,\Gamma}$ be as
above. If $w \in S^\ZD$ contains exactly $c$ forward infinite
disjoint valid paths, then $\htop(X_w,T_{S,\Gamma})=c\log|\Gamma|$.
\end{lem}
\begin{proof}
 With the above notations, since $T$ acts independently on cells in different connected components,
$(X_w,T_{S,\Gamma})$ is isomorphic to $\prod_{i=0}^{c}(X_i,T_i)$,
and so $$\htop(X_w,T_{S,\Gamma})=\sum_{i=0}^{c}\htop(X_i,T_i)$$ By
lemma \ref{lem:no-inifite-path-zero-entropy}, the topological
entropy of
 $(X_{0},T_0)$ is $0$. Each of the other systems $(X_{i},T_i)$ has topological
entropy $\log|\Gamma|$ since it is isomorphic to an inverse limit of
finite extensions of the $\mathbb{N}$-CA on $\Gamma^\NN$ defined by:
$(Tx)_n=x_n+x_{n+1}$.

\end{proof}

In the next section we describe Kari's directed set of tiles from
\cite{kari_xor_snake_1994}, called $S_H$. The tile set $S_H$ is
acyclic and has $0<I(S_H)<\infty$. In fact, we show that this
construction can be carried out in any dimension $d>1$. For this set
of tiles $S_H$ and some fixed group $\Gamma$, denote
$T_H=T_{S_H,\Gamma}$.
 Assuming this, we have our main result:

\begin{thm}\label{thm:main}
For any $d \ge 1$ there exist a surjective $\ZD$-CA with positive,
finite topological entropy.
\end{thm}

\begin{proof}
Let $S_H$ be an acyclic set of directed tiles with
$0<I(S_H)<\infty$.
 By lemma \ref{lem:infinite-path-entropy}, this
means that $\sup_{w \in S_H^{\ZZ^2}}\htop(X_w,T_{H})=
I(S_H)\log|\Gamma| < \infty$. Since $X=\biguplus_{w \in
S_H^{\ZZ^2}}X_w$, we have that $\htop(X,T_H)=\sup_{w \in
S_H^{\ZZ^2}}\htop(X_w,T_{H})$. In section \ref{sec:Hilbet-path} we
will prove the existence of $S_H$ as above, which will complete the
proof.
\end{proof}

Before describing the set of tiles required to complete the proof of
theorem \ref{thm:main}, we not that a slight modification of the
above construction yields the following result about the possible
values the topological entropy of multidimensional CA can obtain:
\begin{prop}\label{prop:CA-entropies-dense}
For any $d \ge 1$, The set of entropies of surjective
$d$-dimensional cellular automata is dense in $[0,\infty)$
\end{prop}
\begin{proof}
Since a product of surjective $d$-dimensional cellular automata also
a surjective $d$-dimensional CA with topological entropy equal to
the sum of the entropies, it is sufficient to prove that there exist
surjective $d$-dimensional CA with arbitrarily small positive
topological entropy.

 Retaining the notations from the beginning of
this section, let $\Gamma$ be some finite group, $S$ an acyclic set
of tiles with $0 < I(S) < \infty$ and $\ZZ_m$ the cyclic group of
order $m$. We define a CA $$T_{S,\Gamma,m}:(S \times \Gamma \times
\ZZ_m) \to (S \times \Gamma \times \ZZ_m)$$ by:
\begin{equation}
\label{eq:kari-ca}
T_{S,\Gamma,m}(x,y,t)_n:=\left\{ \begin{array}{cc} (x_n,y_n+y_{n+d(x_n)},1) & \mbox{if $x$ is valid at $n$ and $t=0$}\\
 (x_n,y_n,t+1) & \mbox{otherwise} \end{array} \right.
 \end{equation}

Evidently, $T_{S,\Gamma,m}$ is surjective (the proof is similar to
lemma \ref{lem:kari-ca-surjective}). Also note that
$T_{S,\Gamma,m}^m$ is conjugate to $T_{S,\Gamma} \times id$, and so
$$\htop(T_{S,\Gamma,m}) =\frac{1}{m}\htop(T_{S,\Gamma})$$
Since $0<\htop(T_{S,\Gamma})<\infty$ and $m$ was arbitrary, the
proof of this proposition is complete.
\end{proof}

\section{\label{sec:Hilbet-path}Kari's tiles and Hilbert space filling paths}
In this section we describe Kari's set of directed tiles $S_H$,
which has the properties announced above. In order to prove our main
result, we only need two properties from this set of tiles $S_H$:
The first property is $0<I(S_H)$, that is, there exist $x \in
S_H^{\ZZ^2}$ with a forward infinite valid path. The second
property, which is harder to prove,  is  $I(S_H)<\infty$. This
property of $S_H$ follows from a lemma proved by Kari (quoted here
as lemma \ref{lem:square-path}). In the first part of this section
(subsection \ref{subsec:hilbert_path}) we describe a substitution
system which is associated with Kari's set of tiles. Subsection
\ref{subsec:kari_tiles} contains a technical description of Kari's
tiles, and it may be skipped by readers who are familiar with
\cite{kari_xor_snake_1994}. The last part of this section contains a
proof that Kari's tiles have a bounded number of forward-infinite
valid paths.

\subsection{\label{subsec:hilbert_path}Hilbert space filling paths}
\begin{figure}
    \input{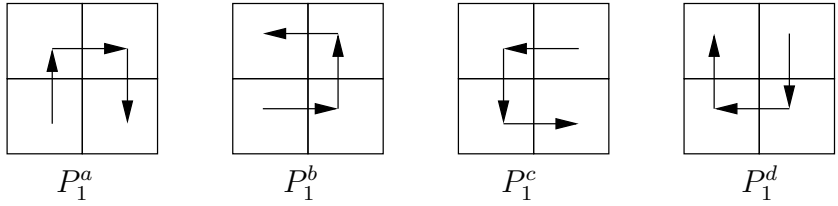}
\caption{The four basic paths through the squares of size $2 \times
2$ tiles.} \label{fig:basic_paths}
\end{figure}

\begin{figure}
    \input{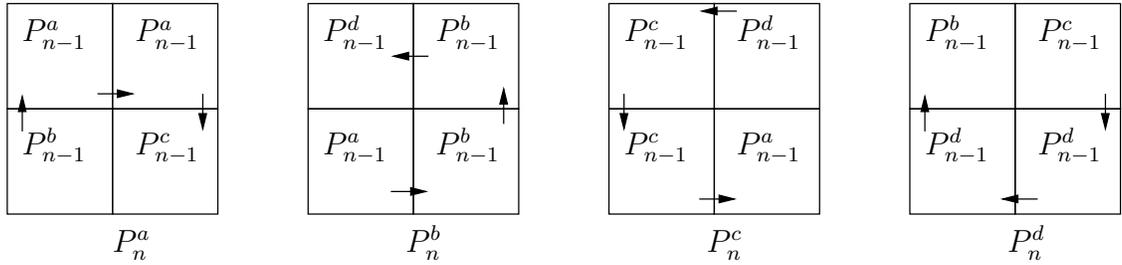}
\caption{The inductive step $n$ in the definition of a Hilbert
path.} \label{fig:inductive_path}
\end{figure}

Before describing Kari's tiles, let us describe a certain
substitution system, which is closely related with these tiles.
 Consider the discrete version of Hilbert's plane
filling curve. This is a path in $\ZZ^2$ which, starting at $(0,0)$
visits each point of the non-negative quarter of $\ZZ^2$ once.
Following is an inductive definition of this path: Define four basic
paths which visit each point in a $2^{n}\times 2^{n}$-square. We
denote these paths by $P_n^a$, $P_n^b$, $P_n^c$, and $P_n^d$.
$P_1^a$, is defined as $(0,0) \rightarrow (0,1) \rightarrow (1,1)
\rightarrow (1,0)$. The other basic paths are obtained by rotations
and reflections. These $4$ basic paths are described in figure
\ref{fig:basic_paths}.
 The path $P_{n+1}^a$ is obtained by walking
according to $P_{n}^b$, $P_{n}^a$, $P_{n}^a$ and $P_{n}^c$.
Similarly, $P_{n+1}^b$, $P_{n+1}^c$ and $P_{n+1}^d$ are defined
inductively as shown in figure \ref{fig:inductive_path}.
\begin{figure}
    \input{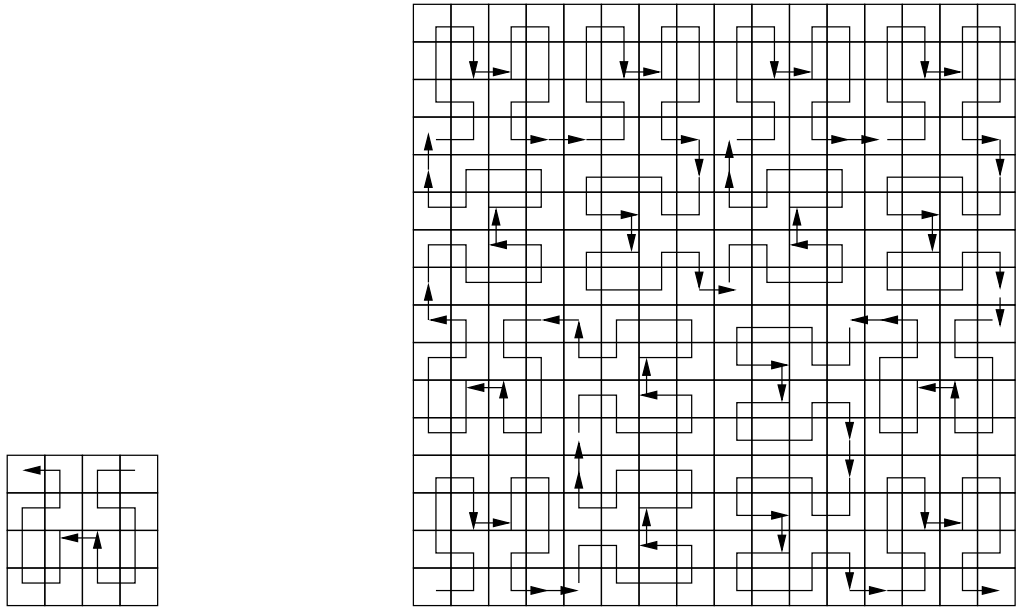}
\caption{Hilbert curves for squares of sizes $2^2\times 2^2$ and
$2^4 \times 2^4$.} \label{fig:hilbert_path}
\end{figure}

A $\ZD$-Hilbert path can also be defined in an analogous manner: the
basic paths are obtained by traversing the vertices of a
$d$-dimensional hypercube using Gray-code. In dimension $d$ there
are $2d$ basic paths, corresponding to the $d$ generators and their
inverses (see Alber and Niedermeier \cite{AN2000} for more on
$\ZD$-Hilbert type paths). For simplicity, we first consider only
the case $d=2$.

\begin{figure}
 \includegraphics[angle=270]{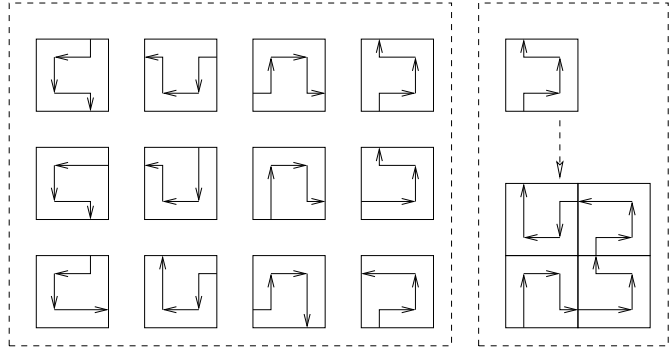}
 \caption{The tiles for the Hilbert substitution system, and a substitution rule.}
 \label{fig:hilbert-tiles}
\end{figure}

We can associate a substitution system to the Hilbert path: there
are $12$ tiles in this system, each corresponding to a basic path
with the inward and outward directions indicated. The substitution
rules reflect the inductive step of the definition of this path. The
tiles and one of the rules of this substitution system is are
described in figure \ref{fig:hilbert-tiles}. The other substitution
rules follow in a symmetric manner, according the definition of the
Hilbert path.
\begin{lem}
The substitution system associated with the Hilbert tiles (described
in figure \ref{fig:hilbert-tiles}) has unique derivation.
\end{lem}
\begin{proof}
It can be directly verified that any tile in the interior of a valid
configuration is a part of a unique path of length $4$ which fills a
$2\times2$ square. Now consider a $\ZZ^2$ configuration in the
associated subshift. The lattice $\ZZ^2$ is partitioned into $2
\times 2$ squares each of which is tiled by a path of length $4$.
This partition uniquely determines the substitution rule applied to
derive the $\ZZ^2$ configuration.
\end{proof}
By the result of Mozes (theorem \ref{thm:subst}), there exist a set
of tiles which implements the Hilbert substitution system. We do not
know however, if any tiling system $S$ which implements the Hilbert
substitution system has the property that $0 < I(S) < \infty$.
Kari's tiles, which are a specific set of tiles which implement this
substitution system, do have this property. We suspect that by
applying the algorithm described by Mozes in \cite{mozes89} on the
Hilbert substitution system, one also ends up with a set of tiles
with the required properties.

\subsection{\label{subsec:kari_tiles}Description of Kari's tiles}
We now describe Kari's tiles, which we denote by $S_H$.  This
directed set of tiles has ``direction'' labels, which are defined in
such a manner that any valid path follows part of the Hilbert path.
The construction is similar to an aperiodic set of tiles constructed
by Robinson \cite{robinson71}, which is also at the heart of Mozes`
result \cite{mozes89}.

\begin{figure}
    \input{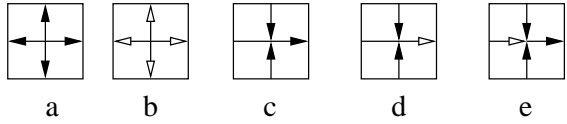}
\caption{The five basic labels of tiles: (a) bold cross (b) blank
cross (c) bold arm (d) blank arm (e) mixed arm.}
\label{fig:basic_tiles}
\end{figure}

Each tile has a \emph{basic label} which is either a \emph{blank
cross}, a \emph{bold cross} a \emph{blank arm}, a \emph{bold arm} or
a \emph{mixed arm}. The five types of basic labels are represented
in figure \ref{fig:basic_tiles}. A tile with a basic label which is
an arm can face in one of the four main directions (north,south,east
or west).
Each cross also has an \emph{orientation label} in NE,NW,SE,SW
refereing to the four cornerwise directions. The tiles labeled with
an arm contain one arrow which faces outwards, and two arrows facing
inwards. We call the arrow facing outwards a \emph{principal arrow}
and the arrows facing inwards \emph{side arrows}. Each arrow in an
arm also contains an orientation label in NE,NW,SE,SW. The principle
arrow of an arm can be labeled with any of the four orientations.

\begin{figure}
    \input{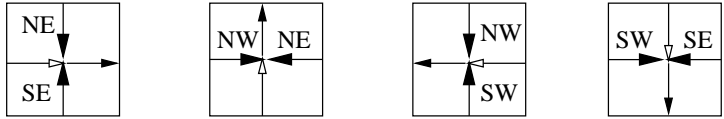}
\caption{The possible orientation labels of the two side arrows of a
mixed arm.}\label{fig:mixed_lables}
\end{figure}

There are restrictions on the orientation labels of the side arrows,
determined by the direction of the arm, and the type of the arm (see
figure \ref{fig:bold_labels}):
\begin{figure}
    \input{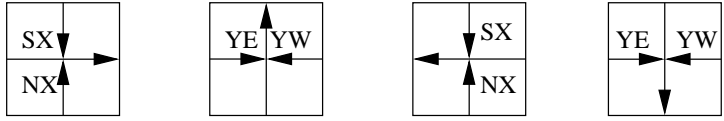}
\caption{The possible orientation labels of the two side arrows of a
bold arm. $X \in \{E,W\}$ and $Y \in \{N,S\}$.}
\label{fig:bold_labels}
\end{figure}
For a bold arm in which the principle arrow is horizontal, there are
two possibilities: either the upper side arrow is labeled with SE
and the lower one with NE, or the upper arrow is labeled with SW and
the lower one with NW. If the arm is vertical, then either the left
arrow is NE and the right arrow is NW or the left arrow is SE and
the right arrow is SW. The possible orientation labels for size
arrows in bold arms are shown in figure \ref{fig:bold_labels}. The
orientation labels for size arrows in a blank arm are similar to a
bold arm (figure \ref{fig:bold_labels}). The possible orientation
labels for arrows in mixed arms are shown in figure
\ref{fig:mixed_lables}.

\begin{figure}
    \input{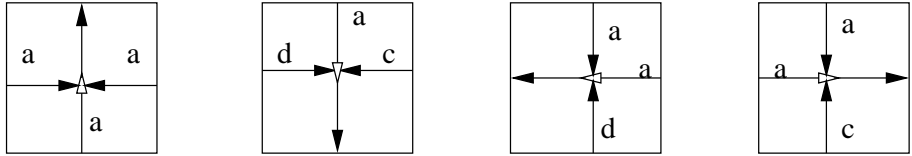}
\caption{The labeling of mixed arms whose principle arrow has the
label $a$.} \label{fig:abcd_mixed_labels}
\end{figure}

Each arrow head also has a \emph{Hilbert-label} in $\{a,b,c,d\}$.
The only restriction imposed on the Hilbert-labels concern crosses
and mixed arms. In each cross all four arrow heads must have the
same Hilbert-label. The Hilbert-labels of the side arrows on the
mixed arms are restricted according to the direction and
Hilbert-label of the principle arrow. The restrictions can be
obtained from figure \ref{fig:inductive_path} as follows: If the
Hilbert-label of the principle arrow is $x \in \{a,b,c,d\}$, and the
principle arrow is facing right, then the upper side arrow has the
label of the path in the upper right corner of the square of
$P^x_{n+1}$, and the lower side arrow has the label of the path in
the lower right corner of the square of $P^x_{n+1}$. For mixed arms
facing the other directions, the restrictions are defined in a
similar manner. For example, figure \ref{fig:abcd_mixed_labels}
 shows the allowed labels of side arrows for mixed arms whose primary
 arrow is labeled $a$.

\begin{figure}
    \input{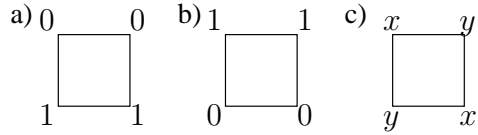}
\caption{The possible corner-parity labels on (a) horizontal arms
(b) vertical arms (d) crosses ($x,y \in \{0,1\}$)
}\label{fig:corner_labels}
\end{figure}

Each tile has four \emph{corner-parity label} in $\{0,1\}$, each
corresponding to a corner of the tile. In a tile which has a
horizontal arm the top corner labels are both $0$, and the bottom
corner labels are both $1$. In a tile which has a vertical arm the
top corner labels are both $1$ and the bottom corner labels are both
$0$. A tile with a cross has two possible labeling for the corner
parities, with opposite corners having the same parity and adjacent
corners having different parities. The possible corner parity labels
are shown in figure \ref{fig:corner_labels}.

There are two more labels for each tile: a \emph{horizontal-parity
label} and a \emph{vertical-parity label}, both of which are in
$\{0,1\}$. The restrictions on these labels are as follows: For a
blank cross, both these parity labels must be $0$. For all other
basic tiles (bold crosses and arms), at least one of these parity
labels must be $1$.

In addition each tile has a direction label, which is in
$\{N,S,E,W\}$. These direction labels determine the direction
function $d:S_H \to \{N,S,E,W\}$, which make $S_H$ a directed set of
tiles.

Let us describe the adjacency rules for these tiles, which define
when a tiling is locally valid. For a configuration to be valid
around a tile the following conditions must hold:
\begin{enumerate}
  \item {\label{rule:arrows_outwards} For any arrow head facing outwards, the arrow head must meet an arrow tail.}
  \item {\label{rule:arrows_inwards} Any arrow head facing inwards must have an arrow head which meets it's tail.}
  \item {\label{rule:labels_match} The orientation labels and Hilbert-labels of meeting arrow heads and tails agree.}
  \item {\label{rule:corner_parity} Corner-parity labels of adjacent corners agree.}
  \item {\label{rule:parity} Both horizontal and vertical parity labels alternate.}
  \item {\label{rule:directions} The direction label is valid, as explained below.}
\end{enumerate}

We now explain the rules which determine when the direction label is
valid (rule \ref{rule:directions}). the direction label of a blank
cross is valid if and only if one of the following conditions are
fulfilled:
\begin{enumerate}
  \item{The direction is N and either the NW neighbor is a bold cross with Hilbert-label $b$ or a vertical arm whose right side arrow has Hilbert-label $a$ or $d$.}
  \item{The direction is W and either the NW neighbor is a bold cross with Hilbert-label $d$ or a horizontal arm whose lower side-arrow has Hilbert-label $c$ or $b$.}
  \item{The direction is S and either the SE neighbor is a bold cross with Hilbert-label $c$ or a vertical arm whose left side arrow has Hilbert-label $a$ or $d$.}
   \item{The direction is W and either the SE neighbor is a bold cross with Hilbert-label $a$ or a horizontal arm whose lower side-arrow has Hilbert-label $c$ or $b$.}
\end{enumerate}
For a tile which is not a blank cross, the direction label is valid
if it agrees with the direction label of a neighboring blank cross
whose direction label points at this tile.

This completes the description of Kari's tiles and their adjacency
rules.

It is possible to extend Kari's construction to a $\ZD$ tiling with
any $d>2$. The details are lengthy but there are no conceptual
difficulties. We comment that these extensions are related to
``multidimensional Hilbert paths''.

\subsection{\label{subsec:finite_infinite_paths}Bounding the number of infinite valid paths}
To complete the proof of theorem \ref{thm:main} for $d=2$, it
remains to show that $0<I(S_H)<\infty$.

\begin{figure}
    \input{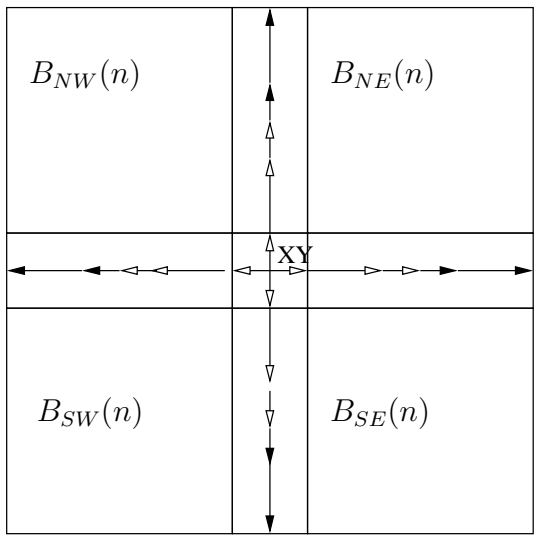}
\caption{Inductive step in definition of $B_{XY}(n+1)$
.}\label{fig:n_square}
\end{figure}

\begin{figure}
    \input{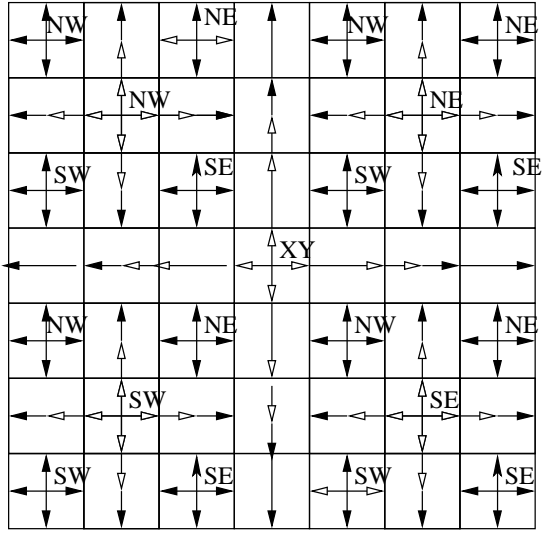}
\caption{$B_{XY}(2)$ with Orientation-labels on the
crosses.}\label{fig:7_square}
\end{figure}

We now show that there exist an infinite valid path for $S_H$. By
induction, we define sequences of valid $S_H$-configurations of
squares of size $(2^{n-1}+1)\times(2^{n-1}+1)$, and denote these by
$B_{XY}(n)$ with $X \in \{N,S\}$ and $Y \in \{E,W\}$. $B_{XY}(0)$ is
a blank cross with orientation-label $XY$. $B_{XY}(n+1)$ is obtained
by surrounding bold cross with orientation label $XY$ by $4$
configurations $B_{SE}(n),B_{SW}(n),B_{NW}(n),B_{NE}(n)$, with arms
labeled correctly as in figure \ref{fig:n_square}. The configuration
$B_{XY}(2)$ with Orientation-labels on the crosses is illustrated in
figure \ref{fig:7_square}. The reader can verify that it is possible
to fill the parity labels and corner-parity labels of these
configurations so they remain valid. This was proved by Kari
\cite{kari_xor_snake_1994}. We now explain how to add Hilbert-labels
and directions to $B_{XY}(n)$ so that the configuration remains
valid. By induction on $n$, we can show that there exists a unique
way to add these labels so that the central cross has Hilbert-label
$x \in \{a,b,c,d\}$. The Hilbert-labels of each of the four
surrounding squares $B_{XY}(n-1)$ are determined by the inductive
definition of the Hilbert path (as in figure
\ref{fig:inductive_path}). The labeling of the paths is determined
by a Hilbert-path which traverses all blank-crosses within this
square. These configurations define arbitrarily long valid paths,
and so  by compactness there exist an infinite valid path for Kari's
tiles (corresponding to a Hilbert path).

Our next goal is to prove that $I(S_H)<\infty$. For this we quote
the following technical lemma, about the structure of valid paths in
the tiling system $S_H$:

\begin{lem} (Kari \cite{kari_xor_snake_1994}, lemma $5$)
\label{lem:square-path} For each $x \in S_H^{\ZZ^2}$ if
$p_1,\ldots,p_N$ is a valid path in $x$ 
, and $N \ge 2\cdot4^{n}$, then there are integers $1 \le i < j \le
N$ with $j-i = 4^n$ so that the path $p_i,\ldots , p_j$ fills a $2^n
\times 2^n$ square.
\end{lem}
 Kari's proof of this lemma involves a delicate examination of
the tiles $S_H$. We comment that it is easier to prove the
corresponding statement about paths which are admissible
configurations for the substitution system $(S_H,\rho_H)$. However,
applying theorem \ref{thm:subst} to deduce the result for the tiling
system is non-trivial, since this lemma refers to configurations of
which only some part (the path) is assumed to be valid. Using lemma
\ref{lem:square-path}, we obtain the following:
\begin{lem}
\label{lem:valid-path-fill-space} Suppose $x \in S_H^{\ZZ^2}$ and
the cells in $p_1,\ldots,p_N$ form a valid path in $x$. There exist
a square  $F \subset \ZZ^s$  centered at $p_1$ which contains all
the cells of the path, such that the path fills up an
$\epsilon$-fraction of the cells in $F$. The constant $\epsilon$ is
independent of $x$ and of the path.
\end{lem}
\begin{proof}
Suppose $p_1, \ldots,p_N$ is a valid path with $2\cdot4^n\le N <
3\cdot4^n$ (if $3\cdot4^n\le N \le 4\cdot4^n$, we look at a prefix).
By lemma \ref{lem:square-path}, there must by some $1\le i_1 \le
4^n$ such that the path $p_{i_1},\ldots,p_{i_1+4^n}$ fills a square
of size $2^{n}\times 2^{n}$. If $i_1>1$ , let $n_1 = \lfloor \log_4
(i_1) \rfloor$. Apply lemma \ref{lem:square-path} on the path
$p_{i_1-4^{n_1}},\ldots,p_{i_1+4^{n-1}}$, and deduce that
$p_{i_1-4^{n_1}},\ldots,p_{i_1}$ fills a square of size
$2^{n_1}\times 2^{n_1}$. Let $i_2=i_1 -4^{n_1}$, continue in this
manner applying lemma \ref{lem:square-path} on path segments of
sizes $2\cdot 4^{n_j}$. Each time set $n_j= \lfloor \log_4 (i_1)
\rfloor$, and  $i_{j+1}=i_j-4^{n_j}$, until $i_j=1$. Each path
segment fills a square of size $2^{n_i}\times 2^{n_i}$. The sizes of
the $n_j$'s are non-increasing, and there can be no more then $3$
$n_j$'s of the same size. Repeat the same procedure on the suffixes
of the path $p_{i_1+4^n},\ldots,p_{N}$, this time the $i_j$'s
increase until $i_j=N$.
 We obtain that $\max_{1\le i,j\le N}
\|p_i-p_j\|_{\infty} \le 6\sum_{i=0}^n 2^i$, and so the entire path
is contained inside a square of size $2^{n+3}\times 2^{n+3}$,
centered at $p_1$, of which it fills at least $2\cdot4^n$ cells.
This completes the proof of the lemma, with $\epsilon =
\frac{4^{n}}{4^{n+3}}=\frac{1}{64}$.
\end{proof}

We now use lemma \ref{lem:valid-path-fill-space} to show that the
number of forward infinite valid paths in any configuration on the
Hilbert tiles is bounded:
\begin{lem}\label{lem:at-most-M-valid-paths}
For each $x \in S_H^{\ZZ^2}$ there exist at most $M= \lceil
\frac{1}{\epsilon} \rceil$ disjoint forward infinite valid paths in
$x$, where $\epsilon$ is the constant from lemma
\ref{lem:valid-path-fill-space}.
\end{lem}
\begin{proof}

Fix some $x \in S_H^{\ZZ^2}$.  Suppose there are more then $M=\lceil
\frac{1}{\epsilon} \rceil$ forward infinite valid paths in $x$, and
fixes cells $c_1,\ldots, c_{M+1}$ with $c_i$ a cell in the $i$'th
path. For any sufficiently large $N$, by lemma
\ref{lem:valid-path-fill-space} each of these paths fill fills up an
$\epsilon$-fraction of the cells in a $2^N \times 2^N$ -square
centered around $c_i$. As $N$ tends to $\infty$, each of these paths
also fill almost an $\epsilon$-fraction of the same square of size
$2^N \times 2^N$. As the paths are disjoint, we reach a
contradiction.
\end{proof}

It follows from lemma \ref{lem:at-most-M-valid-paths} that
$I(S_H)\le M$. On the other hand, we have seen that $I(S_H)>0$. This
concludes the proof if theorem \ref{thm:main}. We remark that it can
actually be proved by a more detailed examination of this system
that there can be at most $4$ disjoint forward infinite valid paths
in any $x \in S_H^{\ZZ^2}$. Also, there exist $x \in S_H^{\ZZ^2}$
with exactly $4$ disjoint infinite valid paths in $x$. Using these
observations one can exactly compute the topological entropy of
$S_H$.

\section{\label{sec:measure-entropy}Measure-theoretic entropy of surjective CA}
This section contains a discussion of surjective CA's as measure
preserving dynamical systems, and the measure-theoretic entropy of
these systems with respect to a ``natural'' measure.

 Denote the symmetric Bernoulli measure by $\mu$: the
state of each cell is distributed uniformly and independently of the
other cells. The following simple proposition was already noted in
the earliest dynamical systems study of CA for the case $\GG=\ZZ$
(see Hedlund's fundamental paper \cite{hedlund69}):
\begin{prop}
If $\GG$ is an amenable group, any surjective $\GG$-CA is measure
preserving with respect to the symmetric Bernoulli measure $\mu$.
\end{prop}

\begin{proof}
First, note that the action $\sigma_g:S^{\GG}\to S^{\GG}$  of $\GG$
on $S^{\GG}$ by translations preserves $\mu$.  The measure theoretic
entropy of this $\GG$-action is equal to the topological entropy of
this action, and $\mu$ is the unique $\sigma$-invariant measure with
this property- $\mu$ is the unique probability measure of maximal
entropy for $\sigma$. Now consider the set $\mathcal{S}=\{\nu \in
\mathcal{P}(S^\GG,\mathcal{B},\sigma):~ \mu=\nu \circ T^{-1}\}$ of
$\sigma$-invariant probability measures. This set $\mathcal{S}$ is
non-empty since $\GG$ is amenable. For any $\nu \in \mathcal{S}$,
$T:(S^{\GG},\mathcal{B},\nu,\sigma) \to
(S^{\GG},\mathcal{B},\mu,\sigma)$ is a measure-theoretic factor map,
and so $h_\nu(S^{\GG},\mathcal{B},\sigma) \ge
h_\mu(S^{\GG},\mathcal{B},\sigma)$, but since $\mu$ is the unique
measure of maximal entropy for $\sigma$, it follows that $\nu=\mu$.
\end{proof}
The closed support of $\mu$ is $S^{\GG}$ for any countable $\GG$,
and so any continuous map which preserves $\mu$ must be surjective.
It follows that for an amenable group $\GG$, a $\GG$-CA is
surjective iff it is $\mu$-preserving.

It is interesting to note that the above characterization of
surjective CA does not hold for general countable groups: Consider
the free group on two generators, denoted by $\mathbb{F}_2$.  We
describe a surjective $\mathbb{F}_2$-CA which does not preserve the
symmetric Bernoulli measure $\mu$ on $\{0,1\}^{\mathbb{F}_2}$:
Denote $2$ generators of $\mathbb{F}_2$ by $a,b$, and consider the
CA $M:\{0,1\}^{\mathbf{F}_2} \to \{0,1\}^{\mathbf{F}_2}$ defined by
the local rule  determined by ``majority vote'' of $x_{wa}$,$x_{wb}$
and $x_{wa^{-1}}$.
\begin{prop}
The the cellular automaton $M:\{0,1\}^{\mathbf{F}_2} \to
\{0,1\}^{\mathbf{F}_2}$ defined as above is surjective, but does not
preserve $\mu$.
\end{prop}

\begin{proof}
To see that $M$ does not preserve $\mu$, let $$A=\{x \in
\{0,1\}^{\mathbf{F}_2} :~ x_{a^{-1}} \ne x_a\}$$

Obviously, $\mu(A)=\frac{1}{2}$. The pre-image of $A$ under $M$ is
all the points $x \in \{0,1\}^{\mathbf{F}_2}$ such that the state of
cell at the identity disagrees either with both the cells at
$a^{-2}$ and $a^{-1}b$ or with both the cells at  $a^2$ and at $ab$:

$$M^{-1}A = \{
x_1\ne x_{a^{-2}}=x_{a^{-1}b}\} \Delta \\
 \{ x_1\ne x_{a^{2}}=x_{ab}\}
$$

It follows that $\mu(A)=\frac{1}{2}$ but $\mu(M^{-1}A)=\frac{1}{8}$,
so $T$ does not preserve the measure $\mu$.

To see that $M$ is surjective, we prove that for any finite $F
\subset \mathbb{F}_2$ and any $y \in \{0,1\}^{\mathbb{F}_2}$ there
exist $x \in \{0,1\}^{\mathbb{F}_2}$ such that $(Tx)_F = y_F$.
Consider the finite graph $G$ with vertex set $F$ and edges
$\{(w_1,w_2) \in F \times F: ~ w_1w_2^{-1} \in
\{a,b,a^{-1},b^{-1}\}\}$. Since $G$ has no cycles, its edges can be
directed so that each connected component is a directed tree. By
adding extra cells to $F$, we can assume that $G$ is a $3$-regular
directed tree. Now we can define $x$: If $w_1,w_2,w_3$ are the
children of $w \in F$, define $x_{w_i}:=y_w$ for $i=1 \ldots 3$.
Define the other states of other cells in $x$ in an arbitrary way.
It follows that indeed $(Tx)_F = y_F$.
\end{proof}

Since the symmetric Bernoulli measure is preserved by any surjective
$\ZD$-CA $T$, one can study the measure theoretic entropy of $T$,
which is bounded by the topological entropy.
\begin{prop} The $\ZZ^2$-CA associated with
Kari's Hilbert-tiles which we denoted by $T_H$ in section
\ref{sec:CA-topo-entropy} has measure theoretic entropy zero with
respect to the symmetric Bernoulli measure.
\end{prop}
\begin{proof}
Recall that for each $w \in S_H^{\ZZ^2}$,  we denoted by $X_w$ all
the points $x \in X$ with directions determined by $w$.
 Since the sets $X_w$ are all $T_H$-invariant, and form a partition of $X$, almost any ergodic component of $\mu$ with respect to
$T_H$ is contained in some $X_w$. Let $F_1$  and $F_2$ be
 squares centered around the origin of dimensions $2^N\times 2^N$
 and $4\cdot2^N \times 4\cdot2^N$ respectively. Suppose there is a valid path staring inside $F_1$ and leaving $F_2$.
   By lemma
\ref{lem:square-path}, such path must fill a square of size $2^N
\times 2^N$ contained inside $F_2$. $F_2$ contains less then
$16\cdot4^N$ squares of size $2^N \times 2^N$, and at least one of
them must be valid. The $\mu$-probability for a square of size
$2^{N} \times 2^{N}$ to be valid is exponentially small in $4^N$,
and so the $\mu$-probability of the event that a valid path starting
in $F_1$ leaves $F_2$ tends to $0$ exponentially as $N \to \infty$.
The Borel-Cantlli lemma implies that almost surely this does not
happen for infinitely many $N$'s, and so with $\mu$ probability $1$
there is no infinite valid path. By lemma
\ref{lem:no-inifite-path-zero-entropy}, when $w$ has no infinite
valid path, $\mu_w$ is supported on a set with zero topological
entropy. We conclude that $\mu$-almost every ergodic component of
$\mu$ has $0$ measure-theoretic entropy.
\end{proof}

In a private communication, Ron Peled suggested that with an
appropriate set of acyclic tiles, $T_{S,\Gamma}$ defined by equation
\eqref{eq:kari-ca} can be a (surjective) CA with finite non-zero
measure theoretic entropy with respect to $\mu $. To see this, apply
the $T_{S,\Gamma}$ construction with the acyclic tile set
$S=\{\uparrow,\rightarrow\}$ with no adjacency restrictions. With
probability $1$ a configuration has $1$ forward infinite path, and
any pair of forward paths eventually coincide. Using lemma
\ref{lem:infinite-path-entropy}  we deduce that $\mu$-almost any
ergodic component of this CA has entropy $\log|\Gamma|$. On the
other hand, this CA has infinite topological entropy.

\section{Concluding remarks and questions}
Our investigation of the existence of a multi-dimensional CA with
finite non-zero entropy was motivated by  the misleading intuition
explained in the introduction. We conclude with some related
questions and remarks.

 Since the CA $T_H$ described in section \ref{sec:CA-topo-entropy} is not injective, a question which follows naturally is:
\begin{question}
 Does there exist an automorphism of a $\ZD$ full shift (an injective $\ZD$-CA) with
 positive, finite topological entropy for $d>1$?
\end{question}

The discussion of measure-theoretic entropy of surjective CA in
section \ref{sec:measure-entropy} raises the following question:
\begin{question}
 Does there exist a $2$-dimensional surjective cellular automaton with positive
measure theoretic entropy (with respect to the symmetric Bernoulli
measure), and finite topological entropy?
\end{question}

  One can study dynamical properties such as entropy  of cellular automata over any countable group.
  We remark that if $\GG$ is a countable group which is not finitely
 generated, $X = \Sigma^{G}$ and $T=X \to X$ is an $\GG$-cellular
 automaton then $h(T) \in \{0,\infty\}$. Here is an explanation of
 this fact: Since $T$ is given by some local rule,
 there exist some finitely generated subgroup $\mathbb{H} < \GG$ such that each cell only interacts via $T$ with other cells
 in the same $\mathbb{H}$-coset. Since $\GG$ is not countably
 generated, $[\GG:\mathbb{H}]=\infty$, and so $(X,T)$ is conjugate
 to the infinite product $\prod_{a \in {\GG /
 \mathbb{H}}}(S^{\mathbb{H}},T)$, and so must have infinite or zero
 entropy. We ask:
\begin{question}
Does there exist a finitely generated, countable group $\GG$ such
that any $\GG$-CA has either zero or infinite topological entropy?
\end{question}

\bibliographystyle{plain}
\bibliography{cellular_entropy}
\end{document}